# Cost – Flow Summation Algorithm Based on Table Form to Solve Minimum Cost – Flow Problem

Eghbal Hosseini[1] [a]

[a]Mechanical and Energy Engineering Department, Erbil Technical Engineering College, Erbil Polytechnic University, Kurdistan Region, Iraq

**Abstract**

The minimum cost – flow problems have been attracted recently in optimization because of their applications in several areas of applied science and real life. Therefore, finding optima solution of these problems would be significant. Although some heuristic approaches have been proposed for solving the problem, but there is no any method to summarize information and converting from a graph form to a table. In this paper, at first all information of the problem are summarized in a table and then an efficient algorithm based on considering costs and flows is proposed. The algorithm is strongly efficient for problems with large size and it has sufficiently suitable results by solving some our generated problems.

**Keywords:** Minimum cost – flow problem; heuristic approaches; random large size.

## 1 Introduction

The minimum cost – flow problem is one of the most important problem in several areas. The objective of this problem is transmit definite products from sources to destinations with minimum cost and satisfaction of the demand and the supply constraints. Two relaxation of this problem is maximal flow and shortest route problems. In other words the problem is a combination of maximal flow and shortest route.

However several heuristic approaches and meta-heuristic algorithms [1-11] have been proposed for solving optimization problems, but there are just few attempts recently to solve minimum cost – flow problems [12-14]. Therefore having a general algorithm to summarize information of the problem and solve it would be remarkable.

Proposing optimal solution needs to start from a feasible solution as an initial basic feasible (IBFS) solution. Therefore initial feasible solution is basic solution which affects to optimal solution for the problem. In other words, finding IBFS would be significant. There are some heuristic methods in references which can find IBFS for the problem [12-14].

[1] Corresponding author
 *E-mail addresses:* eghbal.hosseini@el.epu.edu.iq, kseghbalhosseini@gmail.com

However some heuristic approaches have been proposed in references, but there is no any attempt to prepare general code of algorithms in Matlab. In this paper, the authors will present heuristic table to summarize problem information and also to propose Matlab code to solve large size problems. Therefore many more problems can be solved in the future by these codes.

Most of proposed algorithms have been used to solve small or specific problems. In this paper, using Matlab code several problems in large or small sizes have been generated. Also the algorithm is implemented by few simple MATLAB codes which is based on costs and flows, therefore it has an acceptable computational complexity.

## 2 The Cost – Flow Summation Algorithm

Let a minimum cost – flow problem with n nodes, defined costs and capacities of flows. Make a tableau which consists of the costs and capacities of flows corresponding arcs of the network. In this tableau, each row in the red part relates to capacities of flows from each node to others and the each column in the blue part under the diagonal shows costs of each node to the others. The table proposes the optimal solution if flows in RHS for all rows, except the last one, be zero. General schema of the tableau has been shown in table 1:

Table 1 – General schema of the summarized table of the problem

| 0 | Flow | Flow | RHS |
|---|---|---|---|
| Cost | 0 | Flow | |
| Cost | Cost | 0 | |
| Cost Summation | | | |

### 2.1 Steps of the algorithm

**Step1:**

Make a table includes n rows and n columns which the diagonal is zero. Write costs in 2th row till nth row and flows in 2th column till nth column. If there is no arc between two nodes, for cost write (∞) and for flow write 0. Right Hand Side (RHS) is for flow and in the first row write all flows which should be sent from supply to demand node and RHS of other rows are 0. Below part of the table is for sum of costs.

**Step 2:**

If the solution not be optimal, among rows in RHS except the last row (demand node) choose the largest flow. Corresponding node (rth node) will be selected as sender node.

$$n_r = \max_{1 \le i < n}\{F_i | F_i \text{ is flows in RHS}\}$$

**Step 3:**

Choose the rth column of the cost part in the table and calculate summation of them with costs to the last node one by one. Minimum of these summations will be selected as receiver node using following formula:

$$n_s = \max_{r \le i \le n}\{Sum_i\}$$

$Sum_i$ is summation of cost from $n_r$ to ith node and from it to the demand node.

**Step 4:**

Calculate and rewrite RHS (sum of flows which each node can send). Flow is positive for sender and negative for receivers.

**Step 5:**

Go back to step 2 and repeat the process until finding the optimal solution.

## 3 Computational results

**Example 1:**

Consider following minimum cost – flow problem.

Table 2 – Details of Examples 1

| From Node | To Node | Flow | Cost |
|---|---|---|---|
| 1 | 2 | 7 | 3 |
| 1 | 4 | 5 | 6 |
| 2 | 3 | 7 | 3 |
| 2 | 5 | 3 | 4 |
| 3 | 4 | 6 | 5 |
| 3 | 5 | 4 | 2 |
| 4 | 5 | 8 | 4 |

The best solution of this problem is 12 unit flow by 103$.

Using the proposed algorithm, summation cost – flow table is made:

Step 1:

| Nodes | 1 | 2 | 3 | 4 | 5 | Flow |
|---|---|---|---|---|---|---|
| 1 | 0 | 7 | 0 | 5 | 0 | +12 |
| 2 | 3 | 0 | 7 | 0 | 3 | 0 |
| 3 | ∞ | 3 | 0 | 6 | 4 | 0 |
| 4 | 6 | ∞ | 5 | 0 | 8 | 0 |
| 5 | ∞ | 4 | 2 | 4 | 0 | 0 |
| Cost | | | | | | |

Step 2:

| Nodes | 1 | 2 | 3 | 4 | 5 | Flow |
|---|---|---|---|---|---|---|
| 1 | 0 | 7 | 0 | 5 | 0 | +12 |
| 2 | 3 | 0 | 7 | 0 | 3 | 0 |
| 3 | ∞ | 3 | 0 | 6 | 4 | 0 |
| 4 | 6 | ∞ | 5 | 0 | 8 | 0 |
| 5 | ∞ | 4 | 2 | 4 | 0 | 0 |
| Cost | | 3 | ∞ | 6 | ∞ | 0 |
| Sum of costs | | 7 | ∞ | 10 | ∞ | |

Cost 1= 3*7=21

Step 3:

| Nodes | 1 | 2 | 3 | 4 | 5 | Flow |
|---|---|---|---|---|---|---|
| 1 | 0 | 0 | 0 | 5 | 0 | 0 |
| 2 | 3 | 0 | 7 | 0 | 3 | -7 |
| 3 | ∞ | 3 | 0 | 6 | 4 | 0 |
| 4 | 6 | ∞ | 5 | 0 | 8 | -5 |
| 5 | ∞ | 4 | 2 | 4 | 0 | 0 |
| Cost | | | ∞ | 6 | ∞ | |
| Sum of costs | | | ∞ | 10 | ∞ | |

Cost 2= 6*5=30

Step 4:

| Nodes | 1 | 2 | 3 | 4 | 5 | Flow |
|---|---|---|---|---|---|---|

| Nodes | 1 | 2 | 3 | 4 | 5 | Flow |
|---|---|---|---|---|---|---|
| 1 | 0 | 0 | **0** | 0 | **0** | 0 |
| 2 | 3 | 0 | 7 | **0** | 3 | -4 |
| 3 | ∞ | 3 | 0 | 6 | 4 | 0 |
| 4 | 6 | ∞ | 5 | 0 | 8 | -5 |
| 5 | ∞ | 4 | 2 | 4 | 0 | -3 |
| Cost | | | 3 | ∞ | 4 | |
| Sum of costs | | | 5 | ∞ | **4** | |

Cost 3= 4*3=12

Step 5:

| Nodes | 1 | 2 | 3 | 4 | 5 | Flow |
|---|---|---|---|---|---|---|
| 1 | 0 | 0 | **0** | 0 | **0** | 0 |
| 2 | 3 | 0 | 0 | **0** | 0 | 0 |
| 3 | ∞ | 3 | 0 | 6 | 4 | -4 |
| 4 | 6 | ∞ | 5 | 0 | 8 | -5 |
| 5 | ∞ | 4 | 2 | 4 | 0 | -3 |
| Cost | | | 3 | ∞ | - | |
| Sum of costs | | | 5 | ∞ | | |

Cost 4= 3*4=12

Step 5:

| Nodes | 1 | 2 | 3 | 4 | 5 | Flow |
|---|---|---|---|---|---|---|
| 1 | 0 | 0 | **0** | 0 | **0** | 0 |
| 2 | 3 | 0 | 0 | **0** | 0 | 0 |
| 3 | ∞ | 3 | 0 | 6 | 4 | 0 |
| 4 | 6 | ∞ | 5 | 0 | 8 | -5 |
| 5 | ∞ | 4 | 2 | 4 | 0 | -7 |
| Cost | | | | 5 | 2 | |
| Sum of costs | | | | 9 | **2** | |

Cost 5= 2*4=8

Step 5:

| Nodes | 1 | 2 | 3 | 4 | 5 | Flow |
|---|---|---|---|---|---|---|

| 1 | 0 | 0 | 0 | 0 | 0 | 0 |
|---|---|---|---|---|---|---|
| 2 | 3 | 0 | 0 | 0 | 0 | 0 |
| 3 | ∞ | 3 | 0 | 6 | 0 | 0 |
| 4 | 6 | ∞ | 5 | 0 | 8 | 0 |
| 5 | ∞ | 4 | 2 | 4 | 0 | -12 |
| Cost |  |  |  |  | 4 |  |
| Sum of costs |  |  |  |  | 4 |  |

Cost 6= 4*5=20

Cost =cost1+…+cost6=103

The optimal solution by the proposed algorithm is 103$.

**Example 2:**

Consider following minimum cost – flow problem which has been proposed in DIMACS.

Table 3 – Details of Examples 2

| From Node | To Node | Flow | Cost |
|---|---|---|---|
| 1 | 2 | 4 | 2 |
| 1 | 3 | 2 | 2 |
| 2 | 3 | 2 | 1 |
| 2 | 4 | 3 | 3 |
| 3 | 4 | 5 | 1 |

The best solution using DIMACS is 4 unit flow by 14$.

Using the proposed algorithm, summation cost – flow table is:

Step 1:

| Nodes | 1 | 2 | 3 | 4 | Flow |
|---|---|---|---|---|---|
| 1 | 0 | 4 | 2 | 0 | +4 |
| 2 | 2 | 0 | 2 | 3 | 0 |
| 3 | 2 | 1 | 0 | 5 | 0 |
| 4 | ∞ | 3 | 1 | 0 | 0 |
| Cost |  |  |  |  |  |

Step 2:

| Nodes | 1 | 2 | 3 | 4 | Flow |
|---|---|---|---|---|---|
| 1 | 0 | 4 | 2 | 0 | +2 |
| 2 | 2 | 0 | 2 | 3 | 0 |
| 3 | 2 | 1 | 0 | 5 | -2 |
| 4 | ∞ | 3 | 1 | 0 | 0 |

| | | 2 | 2 | ∞ | |
|---|---|---|---|---|---|
| **Cost** | | | | | |
| **Sum of costs** | | 5 | 3 | ∞ | |

Cost 1= 2*2=4

Step 3:

| **Nodes** | **1** | **2** | **3** | **4** | **Flow** |
|---|---|---|---|---|---|
| 1 | 0 | 4 | 0 | 0 | +2 |
| 2 | 2 | 0 | 2 | 3 | 0 |
| 3 | 2 | 1 | 0 | 5 | -2 |
| 4 | ∞ | 3 | 1 | 0 | 0 |
| **Cost** | | 2 | 2 | ∞ | |
| **Sum of costs** | | | | | |

Cost 2= 2*2=4

Step 4:

| **Nodes** | **1** | **2** | **3** | **4** | **Flow** |
|---|---|---|---|---|---|
| 1 | 0 | 0 | 0 | 0 | 0 |
| 2 | 2 | 0 | 2 | 3 | -2 |
| 3 | 2 | 1 | 0 | 5 | -2 |
| 4 | ∞ | 3 | 1 | 0 | 0 |
| **Cost** | | | 1 | 3 | |
| **Sum of costs** | | | 2 | 3 | |

Cost 3= 1*2=2

Step 5:

| **Nodes** | **1** | **2** | **3** | **4** | **Flow** |
|---|---|---|---|---|---|
| 1 | 0 | 0 | 0 | 0 | 0 |
| 2 | 2 | 0 | 0 | 3 | 0 |
| 3 | 2 | 1 | 0 | 5 | -4 |
| 4 | ∞ | 3 | 1 | 0 | 0 |
| **Cost** | | | | | |
| **Sum of costs** | | | | | |

Cost 4= 1*4=4

Cost =cost1+cost2+cost3+cost4=14

The optimal solution by the proposed algorithm is 14$.

Optimal table:

| Nodes | 1 | 2 | 3 | 4 | Flow |
|---|---|---|---|---|---|
| 1 | 0 | 0 | 0 | 0 | 0 |
| 2 | 2 | 0 | 0 | 3 | 0 |
| 3 | 2 | 1 | 0 | 5 | 0 |
| 4 | ∞ | 3 | 1 | 0 | -4 |
| Cost | | | | | |

**Example 3:**

Consider following minimum cost – flow problem which is a complete graph.

Table 4 – Details of Examples 3

| From Node | To Node | Flow | Cost |
|---|---|---|---|
| 1 | 2 | 15 | 1 |
| 1 | 3 | 5 | 6 |
| 1 | 4 | 30 | 7 |
| 1 | 5 | 40 | 10 |
| 2 | 3 | 10 | 2 |
| 2 | 4 | 15 | 3 |
| 2 | 5 | 15 | 4 |
| 3 | 4 | 7 | 1 |
| 3 | 5 | 10 | 4 |
| 4 | 5 | 20 | 2 |

Some solution to send 20 unit from the origin to the last nodes are: 180$, 200$, 125$ 1nd 120$. Using the proposed algorithm, summation cost – flow table has been shown in the following table also the proposed optimal solution is 120$.

Step 1:

| Nodes | 1 | 2 | 3 | 4 | 5 | Flow |
|---|---|---|---|---|---|---|
| 1 | 0 | 15 | 5 | 30 | 40 | +20 |
| 2 | 1 | 0 | 10 | 15 | 15 | 0 |
| 3 | 6 | 2 | 0 | 7 | 10 | 0 |
| 4 | 7 | 3 | 1 | 0 | 20 | 0 |
| 5 | 10 | 4 | 4 | 2 | 0 | 0 |
| Cost | | | | | | |

To prove feasibility of the proposed algorithm, solving more examples are needed. Therefore, more examples with common size will be solved. Details of Examples 4-9 has been shown in Table 5. Also size and main matrix of the problem, A matrix which includes both costs and flows, have been shown in this Table.

Table 5 – Details of Examples 4-9

| | n | A | | | | | | | | | | | | | | | Cost |
|---|---|---|---|---|---|---|---|---|---|---|---|---|---|---|---|---|---|
| Example 4 | 5 | 0 | 3 | 15 | 15 | 8 | | | | | | | | | | | 172 |
| | | 5 | 0 | 13 | 3 | 7 | | | | | | | | | | | |
| | | 5 | 1 | 0 | 14 | 12 | | | | | | | | | | | |
| | | 5 | 4 | 1 | 0 | 15 | | | | | | | | | | | |
| | | 2 | 3 | 5 | 5 | 0 | | | | | | | | | | | |
| Example 5 | 8 | 0 | 11 | 12 | 5 | 11 | 10 | 3 | 2 | | | | | | | | 187 |
| | | 4 | 0 | 8 | 15 | 6 | 9 | 4 | 12 | | | | | | | | |
| | | 1 | 5 | 0 | 4 | 8 | 11 | 14 | 15 | | | | | | | | |
| | | 5 | 4 | 4 | 0 | 9 | 3 | 3 | 4 | | | | | | | | |
| | | 4 | 2 | 4 | 1 | 0 | 13 | 4 | 13 | | | | | | | | |
| | | 4 | 1 | 2 | 1 | 1 | 0 | 4 | 14 | | | | | | | | |
| | | 5 | 4 | 2 | 5 | 1 | 3 | 0 | 6 | | | | | | | | |
| | | 2 | 4 | 4 | 1 | 3 | 3 | 4 | 0 | | | | | | | | |
| Example 6 | 10 | 0 | 15 | 1 | 12 | 13 | 14 | 2 | 6 | 4 | 13 | | | | | | 224 |
| | | 1 | 0 | 7 | 14 | 3 | 4 | 3 | 3 | 14 | 9 | | | | | | |
| | | 2 | 4 | 0 | 9 | 3 | 13 | 10 | 6 | 8 | 7 | | | | | | |
| | | 3 | 2 | 5 | 0 | 2 | 4 | 2 | 3 | 4 | 7 | | | | | | |
| | | 3 | 3 | 5 | 2 | 0 | 1 | 14 | 15 | 8 | 8 | | | | | | |
| | | 4 | 4 | 2 | 3 | 1 | 0 | 6 | 14 | 6 | 2 | | | | | | |
| | | 1 | 3 | 4 | 5 | 1 | 3 | 0 | 12 | 6 | 4 | | | | | | |
| | | 3 | 1 | 2 | 1 | 4 | 2 | 3 | 0 | 7 | 2 | | | | | | |
| | | 1 | 4 | 2 | 4 | 4 | 4 | 3 | 1 | 0 | 2 | | | | | | |
| | | 2 | 5 | 1 | 5 | 3 | 5 | 1 | 3 | 1 | 0 | | | | | | |
| Example 7 | 15 | 0 | 10 | 8 | 12 | 11 | 14 | 14 | 6 | 11 | 3 | 1 | 12 | 8 | 8 | 14 | 493 |
| | | 5 | 0 | 10 | 10 | 13 | 13 | 9 | 3 | 4 | 14 | 1 | 8 | 3 | 15 | 11 | |
| | | 5 | 3 | 0 | 8 | 8 | 1 | 11 | 1 | 2 | 8 | 2 | 13 | 13 | 11 | 3 | |
| | | 1 | 2 | 2 | 0 | 10 | 8 | 15 | 10 | 13 | 7 | 7 | 13 | 2 | 2 | 3 | |
| | | 5 | 1 | 1 | 1 | 0 | 6 | 13 | 13 | 1 | 6 | 8 | 7 | 10 | 10 | 5 | |
| | | 4 | 4 | 4 | 3 | 3 | 0 | 7 | 1 | 15 | 3 | 2 | 6 | 3 | 8 | 6 | |
| | | 2 | 4 | 1 | 4 | 1 | 2 | 0 | 15 | 14 | 1 | 12 | 5 | 7 | 9 | 15 | |
| | | 4 | 4 | 1 | 5 | 4 | 3 | 3 | 0 | 7 | 15 | 5 | 11 | 10 | 9 | 11 | |
| | | 3 | 2 | 3 | 3 | 5 | 4 | 4 | 2 | 0 | 10 | 3 | 2 | 15 | 3 | 1 | |
| | | 5 | 3 | 2 | 5 | 5 | 3 | 4 | 3 | 2 | 0 | 9 | 14 | 11 | 3 | 6 | |
| | | 2 | 3 | 2 | 5 | 1 | 2 | 1 | 2 | 3 | 2 | 0 | 7 | 15 | 3 | 13 | |
| | | 5 | 3 | 1 | 5 | 5 | 3 | 1 | 2 | 3 | 3 | 2 | 0 | 10 | 6 | 3 | |

|           |    |                                                                 |     |
|-----------|----|-----------------------------------------------------------------|-----|
|           |    | 4 4 2 1 2 2 3 3 1 2 5 1 0 7 8                                   |     |
|           |    | 5 4 3 3 2 3 5 3 3 2 3 4 4 0 2                                   |     |
|           |    | 2 2 5 1 5 5 4 1 2 2 4 1 4 1 0                                   |     |
| Example 8 | 20 | 0 13 12 6 4 12 15 5 11 7 13 12 3 13 15 8 14 9 3 3               | 518 |
|           |    | 3 0 7 12 13 12 5 9 2 2 3 11 8 3 8 3 1 13 9 14                   |     |
|           |    | 2 2 0 11 9 13 14 15 1 13 10 15 8 8 13 4 8 14 9 13               |     |
|           |    | 3 2 2 0 12 9 4 10 2 10 10 11 14 15 12 9 14 9 1 2                |     |
|           |    | 4 2 5 5 0 13 8 13 4 9 10 1 10 6 1 8 3 2 4 3                     |     |
|           |    | 4 2 3 1 5 0 3 1 10 5 9 11 8 9 7 2 8 13 14 5                     |     |
|           |    | 5 5 2 3 1 3 0 4 9 10 7 4 15 2 2 3 3 10 9 1                      |     |
|           |    | 2 1 1 3 1 3 3 0 14 11 12 1 13 15 15 13 12 8 3 6                 |     |
|           |    | 4 4 4 1 1 2 3 4 0 3 1 15 5 5 5 8 10 1 13 9                      |     |
|           |    | 3 5 4 5 3 2 1 4 4 0 13 6 7 1 3 10 5 14 2 15                     |     |
|           |    | 3 1 2 1 2 3 3 3 5 3 0 9 11 15 5 7 7 12 13 2                     |     |
|           |    | 5 4 5 2 4 2 4 4 1 2 2 0 3 6 1 8 6 3 4 14                        |     |
|           |    | 4 5 2 4 1 4 2 5 1 3 3 0 11 8 14 2 12 12 9                       |     |
|           |    | 3 4 2 4 3 1 1 4 3 1 2 4 1 0 3 9 5 3 4 14                        |     |
|           |    | 4 2 5 2 4 1 2 1 3 4 3 3 4 4 0 2 4 1 7 1                         |     |
|           |    | 4 4 5 2 4 2 1 4 3 3 4 4 2 4 3 0 14 3 2 5                        |     |
|           |    | 5 5 2 4 3 3 5 2 2 1 5 4 3 4 3 4 0 7 2 15                        |     |
|           |    | 3 4 3 5 2 1 1 1 3 3 2 4 4 4 5 5 1 0 5 5                         |     |
|           |    | 1 4 1 3 3 5 3 2 4 4 3 2 1 3 2 1 4 2 0 1                         |     |
|           |    | 3 4 2 4 2 4 4 3 1 2 3 2 1 5 3 5 2 4 2 0                         |     |
| Example 9 | 25 | 0 7 1 12 3 3 10 4 5 7 7 6 10 3 3 2 5 12 4 12 11 13 13 5 5       | 655 |
|           |    | 2 0 8 5 13 13 9 4 11 4 7 6 9 15 12 15 4 8 1 12 10 13 15 14 7    |     |
|           |    | 3 5 0 1 9 4 4 5 2 12 12 9 6 13 9 15 14 6 9 6 10 12 12 2 13      |     |
|           |    | 1 1 2 0 1 7 11 12 6 12 14 4 2 4 6 5 14 1 9 3 13 3 8 15 6        |     |
|           |    | 5 2 2 1 0 1 4 6 6 4 15 11 15 7 15 1 10 13 4 14 12 13 9 12 5     |     |
|           |    | 2 2 1 3 1 0 4 5 9 13 5 7 13 10 15 4 13 11 4 8 6 9 13 2 13       |     |
|           |    | 4 1 3 4 2 3 0 13 6 7 9 11 12 12 6 7 15 9 13 5 10 9 15 2 8       |     |
|           |    | 3 1 3 4 2 5 5 0 8 2 14 14 7 12 3 10 4 7 13 3 5 8 6 12 15        |     |
|           |    | 5 3 2 4 2 5 4 4 0 3 4 11 6 15 15 10 13 7 10 15 9 15 11 8 10     |     |
|           |    | 1 1 2 5 5 4 4 2 3 0 14 3 6 15 7 10 14 15 10 2 1 10 9 15 12      |     |
|           |    | 5 3 5 1 2 3 1 4 3 3 0 10 8 4 15 9 1 11 8 1 14 5 4 2 5           |     |
|           |    | 5 5 4 2 3 4 1 1 2 3 5 0 4 10 1 5 5 14 7 12 10 12 2 15 13        |     |
|           |    | 4 2 2 5 5 4 2 1 3 5 3 5 1 0 1 7 5 10 4 9 10 9 7 1 8 7           |     |
|           |    | 3 2 5 1 3 2 5 4 2 2 1 4 3 0 2 7 7 9 13 11 14 1 4 7 15           |     |
|           |    | 1 1 3 5 3 1 1 5 3 2 4 2 3 4 0 12 7 6 1 12 8 3 7 3 12            |     |
|           |    | 5 2 4 5 1 3 3 2 2 4 5 1 4 1 5 0 6 15 1 13 10 9 10 11 2          |     |
|           |    | 5 3 4 3 5 2 1 3 5 1 3 4 1 5 4 3 0 14 1 5 3 14 2 9 10            |     |
|           |    | 1 3 3 5 5 5 4 3 5 3 2 2 4 3 4 5 5 0 10 13 1 13 8 11 4           |     |
|           |    | 3 5 1 1 2 3 3 5 3 3 3 4 1 5 1 3 2 2 0 9 11 15 7 2 1             |     |
|           |    | 4 1 2 1 1 2 5 3 3 1 5 3 5 4 2 5 4 3 5 0 10 12 11 6 15           |     |
|           |    | 5 3 2 3 5 4 2 4 5 3 3 2 1 1 1 3 1 5 4 1 0 8 15 2 4              |     |
|           |    | 1 2 3 1 1 4 2 4 3 3 1 1 1 3 2 5 5 3 2 2 0 12 14 12              |     |
|           |    | 3 4 2 3 4 5 1 4 3 3 5 2 1 4 4 2 5 5 1 2 1 0 5 3                 |     |
|           |    | 4 1 3 5 3 2 2 3 5 1 3 5 1 4 5 2 1 4 5 4 5 1 4 0 13              |     |
|           |    | 2 3 1 5 4 4 5 1 5 3 4 4 5 1 3 4 5 5 5 3 3 5 1 5 0               |     |

To prove efficiency of the proposed algorithm, solving of large size minimum cost – flow problems are required. More large and random size examples have been solved by the algorithm in table 6.

Table 6 – Large and random size problems

|   | Size | Optimal Solution - Cost |
|---|------|-------------------------|
|   |      |                         |

| Example 10 | 50 | 868 |
| Example 11 | 100 | 2282 |
| Example 12 | 300 | 7221 |
| Example 13 | 500 | 11446 |
| Example 14 | 800 | 19476 |
| Example 15 | 1000 | 24640 |

**4 Conclusion and future works**

In this paper, a novel idea has been presented which summarizes all information of the minimum cost – flow problem in a table. Lower computational complexity is because of summarized table and less calculations in the algorithm. Therefore, perhaps this method will be interested by future works in solving much larger problems. Also the algorithm can be suitable tool for solving small problem especially in educational goals. Proposing more efficient algorithms is possible based on the proposed heuristic table in this paper.